\documentclass[10pt,letterpaper,titlepage]{article}
\usepackage{amsmath}
\usepackage{amssymb}
\usepackage[mathscr]{eucal}
\usepackage{verbatim}

\newcommand{\complexes}{\mathbb{C}}

\newcommand{\qed}{\square}

\newcommand{\R}{\mathbb{R}}
\newcommand{\reals}{\R}

\newtheorem{theorem}{Theorem}
\newtheorem{lemma}[theorem]{Lemma}

\newtheorem{corollary}[theorem]{Corollary}
\newtheorem{proposition}[theorem]{Proposition}

\begin{document}
\title
 {$\R^{n} \rtimes G(n)$ is Algebraically Determined}
\author
   {We'am M. Al-Tameemi\\
     Texas A\&M International University\\
     Department of Engineering, Mathematics, and Physics\\
     5201 University Boulevard\\
     Laredo, Texas 78041\\
     Telephone: 956-326-2598\\
     Fax: 956-326-2439\\
     E-mail: weam.altameemi@tamiu.edu\\
     \&\\
     Robert R. Kallman\\
     University of North Texas\\
     Department of Mathematics\\
     1155 Union Circle 311430\\
     Denton, Texas 76203-5017\\
     Telephone: 940-565-3329\\
     Fax: 940-565-4805\\
     E-mail: kallman@unt.edu}

\date\today

\maketitle


\begin{abstract}
\indent
Let $G$ be a Polish (i.e., complete separable metric topological) group.  Define $G$ to be an algebraically
determined Polish group if for any Polish group $L$ and algebraic isomorphism $\varphi: L \mapsto G$, we have
that $\varphi$ is a topological isomorphism.  Let $M(n,\R)$ be the set of $n \times n$ matrices with real
coefficients and let the group $G$ in the above definition be the natural semidirect product $\R^{n} \rtimes
G(n)$, where $n \ge 2$ and $G(n)$ is one of the following groups: either the general linear group $GL(n,\R) =
\left\{ A \in M(n,\R) \ | \ \det(A) \ne 0 \right\}$, or the special linear group $SL(n,\R) = \left\{ A \in
GL(n,\R) \ | \ \det(A) = 1 \right\}$, or $|SL(n,\R)| = \left\{ A \in GL(n,\R) \ | \ |\det(A)| = 1 \right\}$ or
$GL^{+}(n,\R) = \left\{ A \in GL(n,\R) \ | \ \det(A) > 0 \right\}$.  These groups are of fundamental
importance for linear algebra and geometry.  The purpose of this paper is to prove that the natural semidirect
product $\R^{n} \rtimes G(n)$ is an algebraically determined Polish group. Such a result is not true for
$\complexes^{n} \rtimes GL(n,\complexes)$ nor even for $\R^{3} \rtimes SO(3,\R)$.  The proof of this result is
done in a sequence of steps designed to verify the hypotheses of the road map Theorem \ref{theorem_road_map}.
A key intermediate result is that $\varphi^{-1}(SO(n,\R))$ is an analytic subgroup of $L$ for every $n \ge 2$.
\end{abstract}


\newpage

\section{Introduction}

Let $G$ be a Polish (i.e., complete separable metric topological) group.  Define $G$ to be an algebraically
determined Polish group if for any Polish group $L$ and algebraic isomorphism $\varphi: L \mapsto G$, we have
that $\varphi$ is a topological isomorphism. Note that if $G$ is algebraically determined then every algebraic
automorphism of $G$ is continuous.

Let $M(n,\R)$ be the set of $n \times n$ matrices with real coefficients and let the group $G$ in the above
definition be the natural semidirect product $\R^{n} \rtimes G(n)$, where $n \ge 2$ and $G(n)$ is one of the
following groups: either the general linear group $GL(n,\R) = \left\{ A \in M(n,\R) \ | \ \det(A) \ne 0
\right\}$, or the special linear group $SL(n,\R) = \left\{ A \in GL(n,\R) \ | \ \det(A) = 1 \right\}$, or
$|SL(n,\R)| = \left\{ A \in GL(n,\R) \ | \ |\det(A)| = 1 \right\}$ or $GL^{+}(n,\R) = \left\{ A \in GL(n,\R) \
| \ \det(A) > 0 \right\}$.  These groups are of fundamental importance for linear algebra and geometry.  The
purpose of this paper is to prove the following theorem.

\begin{theorem}.\\
Let $L$ be a Polish group and let $\varphi: L \mapsto \R^{n} \rtimes G(n)$ be an algebraic isomorphism. Then
$\varphi$ is a topological isomorphism and therefore $\R^{n} \rtimes G(n)$ is an algebraically determined
Polish group.
\label{theorem_main}
\end{theorem}

The only precedent for this theorem seems to be a result for the $ax + b$ group (\cite{kallman_1984a}).  The
proof of Theorem \ref{theorem_main} of necessity must be rather delicate, for the analogous result is not true
for $(\reals,+)$ ($(\reals,+) \simeq (\reals^{2},+)$) nor for $\complexes^{n} \rtimes GL(n,\complexes)$ (the
field $\complexes$ has $2^{2^{\aleph_{0}}}$ discontinuous automorphisms) nor even for $\R^{3} \rtimes
SO(3,\R)$ (\cite{tits_1972a}).

Throughout this paper we will make free use of basic facts, terminology and notation in descriptive set theory
(\cite{cohn_1980a}, \cite{kechris_1994a}, \cite{kuratowski_1966a}, \cite{moschovakis_1980a},
\cite{parthasarathy_1967a}). In particular, if $X$ is a Polish space then $\mathscr{BP}(X)$ denotes the
$\sigma$-algebra of subsets of $X$ with the Baire property, i.e., the $\sigma$-algebra generated by the Borel
subsets and first category subsets of $X$. Analytic sets are sets with the Baire property.

Kechris and Rosendal \cite{kechris_rosendal_2007a} introduced the notion of Polish groups with ample generics.
This is an important but very special class of Polish groups since any homomorphism of such a Polish group
into a separable topological group is automatically continuous.  Obviously any Polish group with ample
generics is algebraically determined.  The converse is false since it is easy to see that any nontrivial
connected Lie group cannot have ample generics even though there are many examples of Lie groups (e.g., the
real $ax + b$ group) that are algebraically determined.

The proof of Theorem \ref{theorem_main} is done in a sequence of steps following the road map given in the
following theorem.

\begin{theorem} (\cite{jasim_kallman_2014a}, Theorem 4).\\
Let $K$ and $Q$ be two Polish groups and let $\theta: Q \mapsto Aut(K)$ be a group homomorphism that satisfies
$K \times Q \mapsto K$, $(k,q) \mapsto \theta_{q}(k)$ is continuous. Then $K \rtimes_{\theta} Q$ is a Polish
group in the product topology.  Let $L$ be a Polish group and let $\varphi: L \mapsto K \rtimes_{\theta} Q$ be
a group isomorphism. If $\varphi^{-1}(K)$ and $\varphi^{-1}(Q)$ are both analytic subgroups of $L$, then both
$\varphi^{-1}(K)$ and $\varphi^{-1}(Q)$ are closed subgroups of $L$. Next, if, in addition,
$\varphi|_{\varphi^{-1}(K)} : \varphi^{-1}(K) \mapsto K$ is measurable with respect to
$\mathscr{BP}(\varphi^{-1}(K))$, then $\varphi|_{\varphi^{-1}(K)}$ is a topological isomorphism.  Furthermore,
if, in addition, $\theta$ is injective, then $\varphi|_{\varphi^{-1}(Q)} : \varphi^{-1}(Q) \mapsto Q$ is a
topological isomorphism. Finally, under all of these conditions, $\varphi: L \mapsto K \rtimes_{\theta} Q$ is
a topological isomorphism and thus $K \rtimes_{\theta} Q$ is an algebraically determined Polish group.
\label{theorem_road_map}
\end{theorem}

The two steps that are usually the most difficult to verify is that $\varphi^{-1}(Q)$ is an analytic subgroup
of $L$ and that $\varphi|_{\varphi^{-1}(K)} : \varphi^{-1}(K) \mapsto K$ is measurable with respect to
$\mathscr{BP}(\varphi^{-1}(K))$.  We attack these problems in the following sections using the notation of
Theorem \ref{theorem_main}.  Section \ref{section_sln} is devoted to proving that $\varphi^{-1}(SL(n,\R))$ is
an analytic subgroup of $L$ as well noting a few other simple but important observations.


\newpage

\section{$\varphi^{-1}(SL(n,\R))$ is an Analytic Subgroup of $L$}
\label{section_sln}

First note that each natural semidirect product $\R^{n} \rtimes G(n)$ is a Polish topological group.

\begin{lemma}.\\
$\R^{n} \simeq \left\{ (x,I_{n}) \ | \ x \in \R^{n} \right\}$ is a maximal abelian subgroup of $\R^{n} \rtimes
G(n)$. If $L$ is a Polish group and $\varphi: L \mapsto \R^{n} \rtimes G(n)$ is an algebraic isomorphism, then
$\varphi^{-1}(\R^{n})$ is a closed normal maximal abelian subgroup of $L$.
\label{lemma_bag}
\end{lemma}
Proof:\\
\indent
A simple computation shows that $\R^{n}$ is maximal abelian subgroup of $\R^{n} \rtimes G(n)$. Hence,
$\varphi^{-1}(\R^{n})$ is closed subgroup of $L$ since any maximal abelian subgroup of a Hausdorff topological
group is closed.  Therefore $\varphi^{-1}(\R^{n})$ is a closed normal abelian subgroup of $L$.  $\qed$\\

Next notice that $SL(n,\reals) \subseteq G(n) \subseteq \R^{n} \rtimes G(n)$ for each $n \ge 2$.  The rest of this
section is devoted to the proof that $\varphi^{-1}(SL(n,\reals))$ is an analytic subgroup of $L$ for each $n
\ge 2$.  The two cases $n$ is even or $n$ is odd will be discussed separately.

\begin{lemma}.\\
If $G$ is a Polish group and $\emptyset \ne A \subseteq G$ is an analytic set, then $[A,A]$, the subgroup of $G$
generated by commutators of pairs of elements of $A$, is an analytic subgroup of $G$.  If $A$, $B \subseteq G$
are nonempty analytic subsets of $G$, then the subgroup of $G$ algebraically generated by $A$ and $B$ is an
analytic subgroup of $G$.
\label{lemma_box}
\end{lemma}
Proof:\\
\indent
This is elementary using the basic facts that a product of two analytic sets is analytic, that the continuous
image of an analytic set is an analytic set and that the union of a sequence of analytic sets is analytic.
$\qed$

\begin{proposition}.\\
If $-I_{n} \in G(n)$, then the centralizer of $(0,-I_{n})$ in $\R^{n} \rtimes G(n)$ is $\left\{ (0,A) \ | \ A
\in G(n) \right\} \simeq G(n)$.  This is true if $n \ge 2$ is even, if $G(n) = GL(n,\R)$ or if $G(n) =
|SL(n,\R)|$. For such $G(n)$, if $L$ is a Polish group and $\varphi: L \mapsto \R^{n} \rtimes G(n)$ is an
algebraic isomorphism, then $\varphi^{-1} (G(n))$ is a closed subgroup of $L$. In addition under these
circumstances $\varphi^{-1}(SL(n,\R))$ is an analytic subgroup of $L$.
\label{proposition_cow}
\end{proposition}
Proof :\\
\indent
It is an elementary calculation that the centralizer of $(0,-I_{n})$ in $\R^{n} \rtimes G(n)$ is $\left\{ (0,A)  \ | \  A \in G(n) \right\} \simeq G(n)$.
Also,
\begin{eqnarray*}
\varphi^{-1}(G(n)) & = & \varphi^{-1} (\left\{(0,A) \ | \ A \in G(n) \right\}) \\
                   & = & \varphi^{-1} (\text{the centralizer of } \ (0,-I_{n}) \text{ in } \R^{n} \rtimes G(n))\\
                   & = & \text{ the centralizer of } \varphi^{-1} ((0,-I_{n})) \text{ in } L
\end{eqnarray*}
Thus $\varphi^{-1}(G(n))$ is closed in $L$.

Now using the fact that
\begin{eqnarray*}
SL(n,\R) & = & [SL(n,\R),SL(n,\R)]\\
         & = & [GL(n,\R), GL(n,\R)]
\end{eqnarray*}
(\cite{jacobson_1985a}, section 6.7, Lemma 2, page 377) we obtain that
\begin{eqnarray*}
SL(n,\R) & = & [SL(n,\R),SL(n,\R)]\\
         & \subseteq & [G(n),G(n)]\\
         & \subseteq & [GL(n,\R),GL(n,\R)]\\
         & = & SL(n,\R)
\end{eqnarray*}
Therefore $\varphi^{-1}(SL(n,\R)) = \varphi^{-1}([G(n),G(n)]) = [\varphi^{-1}(G(n)),\varphi^{-1}(G(n))]$ is an
analytic subgroup of $L$ by Lemma \ref{lemma_box} since $\varphi^{-1}(G(n))$ is closed in $L$.  $\qed$\\

The analogous result in case $n$ is odd is much more delicate.\\

\begin{lemma}.\\
If $n \ge 3$ is odd, then the centralizer of $\left( 0,\begin{pmatrix} 1 & 0 \\ 0 & -I_{n-1} \end{pmatrix}
\right)$ in $\R^{n} \rtimes G(n)$ is\\
$\mathscr{A} = \left\{ \left( \begin{pmatrix} x_{1} \\ 0 \\ \vdots \\ 0 \end{pmatrix}, \begin{pmatrix} \lambda
& 0 \\ 0 & A \end{pmatrix} \right) \in \R^{n} \rtimes G(n) \right\}$, where $\lambda \ne 0$, $\lambda$ and
$x_{1} \in \R$, $A \in G(n - 1)$ with either $\lambda \cdot \det(A) \ne 0$ if $\R^{n} \rtimes G(n) = \R^{n}
\rtimes GL(n,\R)$, or $\lambda \cdot \det(A) > 0$ if $\R^{n} \rtimes G(n) = \R^{n} \rtimes GL^{+}(n,\R)$, or
$\lambda \cdot \det(A) = 1$ if $\R^{n} \rtimes G(n) = \R^{n} \rtimes SL(n,\R)$ or $|\lambda \cdot \det(A)| =
1$ if $\R^{n} \rtimes G(n) = \R^{n} \rtimes |SL(n,\R)|$.  If $L$ is a Polish group and $\varphi: L \mapsto
\R^{n} \rtimes G(n)$ is an algebraic isomorphism, then $\varphi^{-1}(\mathscr{A})$ is a closed subgroup of
$L$.
\label{lemma_xon}
\end{lemma}
Proof:\\
\indent
If $\left( x, \begin{pmatrix} \lambda & \gamma \\ \beta & A \end{pmatrix} \right) \in \R^{n} \rtimes G(n)$, then
\begin{eqnarray*}
\left( x, \begin{pmatrix} \lambda & \gamma \\ \beta & A \end{pmatrix} \right) \cdot \Big(0, \begin{pmatrix}
1 & 0 \\ 0 & -I_{n-1} \end{pmatrix} \Big) & = & \Big(0,\begin{pmatrix} 1 & 0 \\ 0 & -I_{n-1} \end{pmatrix} \Big)
\cdot \left(x, \begin{pmatrix} \lambda & \gamma \\ \beta & A \end{pmatrix} \right)\\
\text{if and only if } \left(x, \begin{pmatrix} \lambda & \gamma \\ \beta & A \end{pmatrix} \cdot \begin{pmatrix}
1 & 0 \\ 0 & -I_{n - 1} \end{pmatrix} \right) & = & \left( \begin{pmatrix} 1 & 0 \\ 0 & -I_{n-1} \end{pmatrix} (x),
\begin{pmatrix} 1 & 0 \\ 0 & -I_{n - 1} \end{pmatrix} \cdot \begin{pmatrix} \lambda & \gamma \\ \beta & A
\end{pmatrix}  \right)\\
\text{if and only if } \left( \begin{pmatrix} x_{1} \\ x_{2} \\ \vdots \\ x_{n} \end{pmatrix}, \begin{pmatrix}
\lambda & -\gamma \\ \beta & -A \end{pmatrix} \right) & = & \left( \begin{pmatrix} x_{1} \\ -x_{2} \\ \vdots \\ -x_{n}
\end{pmatrix},\begin{pmatrix} \lambda & \gamma \\ -\beta & -A \end{pmatrix} \right).
\end{eqnarray*}

\noindent
Therefore $\beta = 0$, $\gamma = 0$ and $x_{i}= 0$ for $ 2 \le i \le n$ and
$\left( \begin{pmatrix} x_{1} \\ x_{2} \\ \vdots \\ x_{n} \end{pmatrix}, \begin{pmatrix} \lambda &
\gamma \\ \beta & A \end{pmatrix} \right) = \left( \begin{pmatrix} x_{1} \\ 0 \\ \vdots \\ 0 \end{pmatrix},
\begin{pmatrix} \lambda & 0 \\ 0 & A \end{pmatrix} \right)$\\
\noindent
with $\lambda$, $x_{1} \in \R$, $A \in GL(n - 1)$ and $\lambda \cdot \det(A)$ is any one of the four mentioned
cases.\\

Now, $\varphi^{-1}(\mathscr{A}) = \varphi^{-1} (\text{the centralizer of } \left( 0, \begin{pmatrix} 1 & 0 \\
0 & -I_{n - 1} \end{pmatrix} \right) \text{ in } \R^{n} \rtimes G(n)) = \text{the centralizer of }$\\ $
\varphi^{-1}(( 0, \begin{pmatrix} 1 & 0 \\ 0 & -I_{n-1} \end{pmatrix}))$ in $L$, which is a closed subgroup of
$L$. $\qed$\\

\begin{lemma}.\\
If $n \ge 3$ is odd and $\mathscr{A}$ is as in Lemma \ref{lemma_xon}, then in all four cases the commutator
subgroup of $\mathscr{A}$ is\\
\noindent
$\mathscr{D} = \left\{ \left( \begin{pmatrix} x_{1} \\ 0 \\ \vdots \\ 0 \end{pmatrix}, \begin{pmatrix} 1 & 0 \\
0 & D \end{pmatrix} \right) \ | \ x_{1} \in \R \text{, } D \in SL(n - 1,\R) \right\}$.
If $L$ is a Polish group and $\varphi: L \mapsto \R^{n} \rtimes G(n)$ is an algebraic isomorphism, then
$\varphi^{-1} (\mathscr{D})$ is an analytic subgroup of $L$.
\label{lemma_zon}
\end{lemma}
Proof:\\
\indent
Compute that
$\left[ \left( \begin{pmatrix} x' \\ 0 \\ \vdots \\ 0 \end{pmatrix}, \begin{pmatrix} \lambda' & 0 \\ 0 & A'
\end{pmatrix} \Big ), \Big( \begin{pmatrix} x \\ 0 \\ \vdots \\ 0 \end{pmatrix}, \begin{pmatrix}
\lambda & 0 \\ 0 & A \end{pmatrix} \right) \right]$ =
$\left( \begin{pmatrix} x_{1} \\ 0 \\ \vdots \\ 0 \end{pmatrix}, \begin{pmatrix} 1 & 0 \\ 0 & D
\end{pmatrix} \right)$,\\
where $x_{1} = (-1 + \lambda') x + (1 - \lambda) x'$ and $A'AA'^{-1}A^{-1} = D \in SL(n - 1,\reals)$.\\

On the other hand elements of the form $\left( \begin{pmatrix} x_{1} \\ 0 \\ \vdots \\ 0 \end{pmatrix},
\begin{pmatrix} 1 & 0 \\ 0 & D \end{pmatrix} \right)$ are in the commutator subgroup of $\mathscr{A}$ since
if we take $\lambda' = 1/2$, $\lambda = 1/2$, $x = -x_{1}$, $x' = x_{1}$ and $A$ and $A'$ suitable simple
diagonal matrices (remember $n \ge 3$), then $\left(\begin{pmatrix}  x_{1} \\ 0 \\ \vdots \\ 0
\end{pmatrix},\begin{pmatrix} 1 & 0 \\ 0 & I_{n-1} \end{pmatrix} \right)$ is in the commutator subgroup of
$\mathscr{A}$. Also, $\left(\begin{pmatrix} 0 \\ 0 \\ \vdots \\ 0 \end{pmatrix},\begin{pmatrix} 1 & 0 \\ 0 & D
\end{pmatrix} \right)$, where $D \in SL(n - 1,\R)$, is in the commutator subgroup of $\mathscr{A}$
(\cite{jacobson_1985a}).  Therefore, the product\\
$\left(\begin{pmatrix}  x_{1} \\ 0 \\ \vdots \\ 0 \end{pmatrix}, \begin{pmatrix} 1 &
0 \\ 0 & I_{n-1} \end{pmatrix} \right) \cdot \left( \begin{pmatrix} 0 \\ 0 \\ \vdots \\ 0
\end{pmatrix},\begin{pmatrix} 1 & 0 \\ 0 & D \end{pmatrix} \right) = \left( \begin{pmatrix} x_{1} \\ 0 \\
\vdots \\ 0 \end{pmatrix},\begin{pmatrix} 1 & 0 \\ 0 & D \end{pmatrix} \right)$ is in the commutator subgroup
of $\mathscr{A}$.

Finally, $\varphi^{-1}(\mathscr{D}) = \varphi^{-1}( \text{ the commutator subgroup of } \mathscr{A} ) =
\text{ the commutator subgroup of } \varphi^{-1}(\mathscr{A})$, which is an analytic subgroup of $L$ by Lemma
\ref{lemma_box}.   $\qed$

\begin{lemma}.\\
If $n \ge 3$ is odd, then the commutator subgroup of
$\mathscr{D} = \left\{\left( \begin{pmatrix} x_{1} \\ 0 \\ \vdots \\ 0 \end{pmatrix}, \begin{pmatrix}
1 & 0 \\ 0 & D \end{pmatrix} \right) \ | \ x_{1} \in \R \text{, } D \in SL(n - 1,\R) \right\}$\\
is
$\mathscr{B} = \left\{\left( \begin{pmatrix} 0 \\ 0 \\ \vdots \\ 0 \end{pmatrix}, \begin{pmatrix} 1 & 0 \\ 0 &
B \end{pmatrix} \right) \ | \ B \in SL(n - 1,\R) \right\}$ and $\varphi^{-1}(\mathscr{B})$ is an analytic
subgroup of $L$.\\
Also the commutator subgroup of $\mathscr{D'} = \left\{ \left( \begin{pmatrix} 0 \\ 0 \\
\vdots \\ x_{n} \end{pmatrix}, \begin{pmatrix} D & 0 \\ 0 & 1 \end{pmatrix} \right) \ | \ x_{n} \in \R
\text{, } D \in SL(n - 1,\R) \right\}$ is\\
\noindent
$\mathscr{B'} = \left\{\left( \begin{pmatrix} 0 \\ 0 \\ \vdots \\ 0 \end{pmatrix}, \begin{pmatrix} B & 0 \\ 0
& 1 \end{pmatrix} \right) \ | \ B \in SL(n - 1,\R) \right\}$ and $\varphi^{-1}(\mathscr{B'})$ is an analytic
subgroup of $L$.
\label{lemma_nur}
\end{lemma}
Proof:\\
\indent
Compute that $\left[ \left( \begin{pmatrix} x_{1} \\ 0 \\ \vdots \\ 0 \end{pmatrix}, \begin{pmatrix} 1 & 0 \\
0 & D \end{pmatrix} \right), \left( \begin{pmatrix}  x'_{1} \\ 0 \\ \vdots \\ 0 \end{pmatrix}, \begin{pmatrix}
1 & 0 \\ 0 & D' \end{pmatrix} \right) \right] = \left( \begin{pmatrix} 0 \\ 0 \\ \vdots \\ 0 \end{pmatrix},
\begin{pmatrix} 1 & 0 \\ 0 & B \end{pmatrix} \right)$, where \\ $B = DD'D^{-1} D'^{-1} \in SL(n - 1,\R)$. On
the other hand elements of the form $\left( \begin{pmatrix} 0 \\ 0 \\ \vdots \\ 0 \end{pmatrix},
\begin{pmatrix} 1 & 0 \\ 0 & B \end{pmatrix} \right)$ are in the commutator subgroup of $\mathscr{D}$ by
\cite{jacobson_1985a}.  Now $\varphi^{-1}(\mathscr{B}) = \varphi^{-1} \left(\text{the commutator subgroup of }
\mathscr{D} \right) =\\
\text{the commutator subgroup of } \varphi^{-1}(\mathscr{D})$, which is an analytic subgroup of $L$ by Lemma
\ref{lemma_box}.

Similarly one can show that the commutator subgroup of $\mathscr{D'}$ is $\mathscr{B'}$ and therefore
$\varphi^{-1}(\mathscr{B'})$  is an analytic subgroup of $L$.  $\qed$

\begin{lemma}.\\
If $n \ge 3$ is odd, let $\mathscr{B}$ and $\mathscr{B'}$ be as in Lemma \ref{lemma_nur} and let $\mathscr{K}$
be the group generated by $\mathscr{B}$ and $\mathscr{B'}$.  Then $\varphi^{-1}(\mathscr{K})$ is an analytic
subgroup of $L$.
\label{lemma_kit}
\end{lemma}
Proof:\\
\indent
Both $\varphi^{-1}(\mathscr{B})$ and $\varphi^{-1}(\mathscr{B'})$ are analytic subgroups of $L$ by Lemma
\ref{lemma_nur}.  Therefore, the group generated by them, namely $\varphi^{-1}(\mathscr{K})$, is also an
analytic subgroup of $L$ by Lemma \ref{lemma_box}.   $\qed$

\begin{lemma}.\\
If $n \ge 2$ view $\R^{n} - \left\{ 0 \right\}$ as an $SL(n,\R)$-space and let $z = (1, 0, \cdots, 0) \in
\R^{n} - \left\{ 0 \right\}$. Then the stability group $SL(n,\R)_{z}= \left\{ a \in
SL(n,\R) \ | \ a \cdot z = z \right\}$ is contained in the group generated by elements of the form\\
$\left\{ \begin{pmatrix} 1 & 0 \\ 0 & B \end{pmatrix} \mid \ B \in SL(n - 1,\R) \right\}$ and $\left\{
\begin{pmatrix} B & 0 \\ 0 & 1 \end{pmatrix} \mid \ B \in SL(n - 1,\R) \right\}$.
\label{lemma_hey}
\end{lemma}
Proof:\\
\indent
If $\lambda_{1}$, $\lambda_{2} \in \reals$, compute that\\
$\begin{pmatrix} 1 & 0 & \cdots & \lambda_{1} & 0 \\ 0 & 1 & \cdots & 0 & 0 \\ \vdots & \vdots & \ddots &
\vdots & \vdots \\ 0 & 0 & \cdots & 1 & 0 \\ 0 & 0 & \cdots & 0 & 1 \end{pmatrix} \cdot \begin{pmatrix} 1 & 0
& \cdots & 0 & 0 \\ 0 & 1 & \cdots & 0 & 0 \\ \vdots & \vdots & \ddots & \vdots & \vdots \\ 0 & 0 & \cdots & 1
& \lambda_{2} \\ 0 & 0 & \cdots & 0 & 1 \end{pmatrix}  = \begin{pmatrix} 1 & 0 & \cdots & \lambda_{1} &
\lambda_{1} \lambda_{2} \\ 0 & 1 & \cdots & 0 & 0 \\ \vdots & \vdots & \ddots & \vdots & \vdots \\ 0 & 0 &
\cdots & 1 & \lambda_{2} \\ 0 & 0 & \cdots & 0 & 1 \end{pmatrix}$ and\\
$\begin{pmatrix} 1 & 0 & \cdots & \lambda_{1} & \lambda_{1} \lambda_{2} \\ 0 & 1 & \cdots & 0 & 0 \\ \vdots &
\vdots & \ddots & \vdots & \vdots \\ 0 & 0 & \cdots & 1 & \lambda_{2} \\ 0 & 0 & \cdots & 0 & 1 \end{pmatrix}
\cdot \begin{pmatrix} 1 & 0 & \cdots & - \lambda_{1} & \lambda_{1} \lambda_{2} \\ 0 & 1 & \cdots & 0 & 0 \\
\vdots & \vdots & \ddots & \vdots & \vdots \\ 0 & 0 & \cdots & 1 & - \lambda_{2} \\ 0 & 0 & \cdots & 0 & 1
\end{pmatrix} = \begin{pmatrix} 1 & 0 & \cdots & 0 & \lambda_{1} \lambda_{2} \\ 0 & 1 & \cdots & 0 & 0 \\
\vdots & \vdots & \ddots & \vdots & \vdots \\ 0 & 0 & \cdots & 1 & 0 \\ 0 & 0 & \cdots & 0 & 1
\end{pmatrix}$.

An element $a \in SL(n,\R)_{z}$ is of the form  $\begin{pmatrix} 1 & a_{12} & \cdots & a_{1(n-1)}
& a_{1n} \\ 0 & \ & \  & \  & \  \\ \vdots & \  & A & \  & \  \\ 0 & \  & \  & \  & \  \end{pmatrix}$,
where $A \in SL(n - 1,\R)$.\\
Observe that
$\begin{pmatrix} 1 & a_{12} & \cdots & a_{1(n-1)} & 0 \\ 0 & 1 & \cdots & 0 & 0 \\ \vdots & \vdots & \ddots &
\vdots & \vdots \\ 0 & 0 & \cdots & 1 & 0 \\ 0 & 0 & \cdots & 0 & 1 \end{pmatrix} \cdot \begin{pmatrix} 1 & 0
& \cdots & 0 & \lambda_{1} \lambda_{2} \\ 0 & 1 & \cdots & 0 & 0 \\ \vdots & \vdots & \ddots & \vdots & \vdots
\\ 0 & 0 & \cdots & 1 & 0 \\ 0 & 0 & \cdots & 0 & 1 \end{pmatrix} = \begin{pmatrix} 1 & a_{12} & \cdots &
a_{1(n-1)} & \lambda_{1} \lambda_{2} \\ 0 & 1 & \cdots & 0 & 0 \\ \vdots & \vdots & \ddots & \vdots & \vdots
\\ 0 & 0 & \cdots & 1 & 0 \\ 0 & 0 & \cdots & 0 & 1 \end{pmatrix}$.
Hence $\begin{pmatrix} 1 & 0 \\ 0 & A \end{pmatrix} \cdot  \begin{pmatrix} 1 & a_{12} & \cdots & a_{1(n-1)} &
\lambda_{1} \lambda_{2} \\ 0 & \ & \  & \  & \  \\ \vdots & \  & I_{n-1} & \  & \  \\ 0 & \  & \  & \  & \
\end{pmatrix} = \begin{pmatrix} 1 & a_{12} & \cdots & a_{1(n-1)} & \lambda_{1} \lambda_{2} \\ 0 & \ & \  & \ &
\  \\ \vdots & \  & A & \  & \  \\ 0 & \  & \  & \  & \  \end{pmatrix}$, where\\
$A \in SL(n - 1,\R)$.  So every element of the stability subgroup is in the subgroup generated by elements of
the form $\left\{ \begin{pmatrix} 1 & 0 \\ 0 & B \end{pmatrix} \ | \ B \in SL(n - 1,\R) \right\}$ and $\left\{
\begin{pmatrix} B & 0 \\ 0 & 1 \end{pmatrix} \ | \ B \in SL(n - 1,\R) \right\}$.  $\qed$

\begin{lemma}.\\
If $n \ge 2$, then $SL(n,\R)$ acts transitively on $\R^{n} - \left\{ 0 \right\}$ and $SL(n,\R)/SL(n,\R)_{z}
\simeq  \R^{n} - \left\{ 0 \right\}$.
\label{lemma_how}
\end{lemma}
Proof:\\
\indent
This is well known and is an easy exercise.  $\qed$

\begin{lemma}.\\
If $n \ge 3$ is odd, let $\mathscr{K}$ be the subgroup of $SL(n,\R)$ generated by
$\left\{ \begin{pmatrix} 1 & 0 \\ 0 & B \end{pmatrix} \ | \ B \in SL(n - 1,\R) \right\}$ and $\left\{
\begin{pmatrix} B & 0 \\ 0 & 1 \end{pmatrix} \ | \ B \in SL(n - 1,\R) \right\}$.
Then $\mathscr{K}$ acts transitively on $\R^{n} - \left\{ 0 \right\}$.
\label{lemma_who}
\end{lemma}
Proof:\\
\indent
Let $x = (x_{1}, x_{2},\cdots, x_{n}) \ne 0$ and $z = (1,0, \cdots,0)$. It is suffices to show there exists $g
\in \mathscr{K}$ such that $g \cdot x = z$. Use Lemma \ref{lemma_how}.  There exists $\begin{pmatrix} 1 & 0 \\
0 & B \end{pmatrix} \in SL(n,\R)$, where $B \in SL(n - 1,\R)$, such that $\begin{pmatrix} 1 & 0 \\ 0 & B
\end{pmatrix} \cdot \begin{pmatrix} x_{1} \\ x_{2} \\ x_{3} \\ \vdots \\ x_{n} \end{pmatrix} = \begin{pmatrix}
x_{1} \\ x'_{2} \\ 0 \\ \vdots \\ 0 \end{pmatrix}$.  Also $\begin{pmatrix} B' & 0 \\ 0 & 1 \end{pmatrix} \cdot
\begin{pmatrix} x_{1} \\ x'_{2} \\ 0 \\ \vdots \\ 0 \end{pmatrix} = \begin{pmatrix} 1 \\ 0 \\ 0 \\ \vdots \\ 0
\end{pmatrix}$ for some $B' \in SL(n - 1,\reals)$. Therefore $\mathscr{K}$ acts transitively on $\R^{n} -
\left\{ 0 \right\}$.  $\qed$

\begin{lemma}.\\
If $G$ is a group and $H \subset K \subset G$ are subgroups such that $K$ acts transitively on $G/H$, then $K=
G$.
\label{lemma_ruf}
\end{lemma}
Proof:\\
\indent
If $a \in G$, then $a \in aH$.  Since $K$ acts transitively on $G/H$, there exists $k \in K$ such that $a \in
aH = k \cdot H \subset KH \subset K \cdot K = K$. Therefore $K = G$.  $\qed$

\begin{proposition}.\\
If $n \ge 3$ is odd, then $SL(n,\R)$ is the group generated by elements of the form
$\left\{ \begin{pmatrix} 1 & 0 \\ 0 & B \end{pmatrix} \mid \ \ B \in SL(n-1, \R) \right\}$ and $\left\{
\begin{pmatrix} B & 0 \\ 0 & 1 \end{pmatrix} \mid \ B \in SL(n-1, \R) \right\}$.
If $L$ is a Polish group and $\varphi: L \mapsto \R^{n} \rtimes G(n)$ is an algebraic
isomorphism, then $\varphi^{-1}(SL(n,\R))$ is an analytic subgroup of $L$.
\label{proposition_dar}
\end{proposition}
Proof:\\
\indent
Let $\mathscr{K}$ be the subgroup olf $SL(n,\reals)$ generated by elements of the form $\left\{
\begin{pmatrix} 1 & 0 \\ 0 & B \end{pmatrix} \mid \ B \in SL(n-1, \R) \right\}$ and $\left\{ \begin{pmatrix} B
& 0 \\ 0 & 1 \end{pmatrix} \ | \ B \in SL(n - 1,\R) \right\}$. Then $SL(n,\R)_{z} \subset \mathscr{K}$ by
Lemma \ref{lemma_hey} and $\mathscr{K}$ acts transitively on\\ $SL(n,\R)/SL(n,\R)_{z}$ by Lemma
\ref{lemma_who}. $SL(n,\R) = \mathscr{K}$ by Lemma \ref{lemma_ruf}. Therefore $\varphi^{-1}(SL(n,\R))=
\varphi^{-1}(\mathscr{K})$ is an analytic subgroup of $L$ by Lemma \ref{lemma_kit}.  $\qed$


\newpage

\section{$\varphi^{-1}(G(n))$ is an Analytic Subgroup of $L$}
\label{section_gn}

The goal of this section is to prove that if $L$ is a Polish group and $\varphi: L \mapsto \R^{n} \rtimes
G(n)$ is an algebraic isomorphism, then $\varphi^{-1}(G(n))$ is an analytic subgroup of $L$.  This has already
been done in the cases $G(n) = GL(n,\R)$ ($n \ge 2$ even), $G(n) = |SL(n,\R)|$ ($n \ge 2$), $G(n) = SL(n,\R)$
($n \ge 2$) and $G(n) = GL^{+}(n,\R)$ ($n \ge 2$ even).  The only case remaining is $G(n) = GL^{+}(n,\R)$ and
$n \ge 3$ is odd.  The results of the previous section on $SL(n,\R)$ will lead to a quick proof.  The
following lemma certainly is well known.

\begin{lemma}.\\
The centralizer of $SL(n,\R)$ in $GL(n,\R)$ is the group of nonzero scalar matrices.
\label{lemma_swing}
\end{lemma}
Proof:\\
\indent
A scalar matrix will commute with any matrix in $GL(n,\R)$. Conversely, let $A = (a_{ij})$ belong to the
centralizer of $SL(n,\R)$ in $GL(n,\R)$. Write $E_{ij}$ for the elementary $n \times n$ matrix with $1$ in the
$(i,j)$th position and $0$ elsewhere. Now $I + E_{ij}$ in $SL(n,\R)$ if $i \ne j$ and so $A$ and $E_{ij}$
commute.  The $(k,j)$th entry of $AE_{ij}$ is $a_{ki}$ while that of $E_{ij}A$ is $0$ if $k \ne i$ and
is $a_{jj}$ if $k = i$.  Hence $a_{ki} = 0$ if $k \ne i$ and $a_{ii} = a_{jj}$, which shows that $A$ is a
scalar.    $\qed$

\begin{lemma}.\\
For $n \ge 2$ the centralizer of $SL(n,\R)$ in $G(n)$, say $C$, is the set of diagonal elements of $G(n)$.  If
$G(n) = GL^{+}(n,\R)$, then $GL^{+}(n,\R) = C \cdot SL(n,\R)$.
\end{lemma}
Proof:\\
\indent
The first assertion folows from Lemma \ref{lemma_swing}. For the second, suppose $T \in GL^{+}(n,\R)$ and
$\det(T) = \lambda^{n} > 0$.  Let $D(\lambda)$ be the diagonal matrix with $\lambda$ as the diagonal entries.
$D(\lambda) \in C$ and $T = D(\lambda) \cdot D(\lambda)^{-1}T \in C \cdot SL(n,\R)$.  $\qed$

\begin{proposition}.\\
Let $L$ be a Polish group and $\varphi: L \to \R^{n} \rtimes GL^{+}(n,\R)$ be an algebraic isomorphism.  Then
$\varphi^{-1}(GL^{+}(n,\R))$ is an analytic subgroup of $L$.
\label{proposition_rat}
\end{proposition}
Proof:\\
\indent
Let $C$ be the centralizer of $SL(n,\R)$ in $\R^{n} \rtimes GL^{+}(n,\R)$.  Then $C \subset GL^{+}(n,\R)$,
$GL^{+}(n,\R) = C \cdot SL(n,\R)$ and $\varphi^{-1}(C)$ is closed in $L$ since it is the centralizer of
$\varphi^{-1}(SL(n,\R))$ in $L$.  Hence, $\varphi^{-1}(GL^{+}(n,\R)) = \varphi^{-1}(C \cdot SL(n,\R)) =
\varphi^{-1}(C) \cdot \varphi^{-1}(SL(n,\R))$ is a product of two analytic subsets of $L$ and therefore itself
is an analytic subset of $L$.  $\qed$

\begin{corollary}.\\
If $L$ is a Polish group and $\varphi: L \mapsto \R^{n} \rtimes G(n)$ is an algebraic isomorphism, then
$\varphi^{-1}(G(n))$ is an analytic subgroup of $L$.
\label{corollary_gn}
\end{corollary}
Proof:\\
\indent
This follows from the comments and results of this section.  $\qed$


\newpage

\section{$\varphi^{-1}(SO(n,\R))$ is an Analytic Subgroup of $L$}
\label{section_so}

The goal of this section is to prove that if $L$ is a Polish group and $\varphi: L \mapsto \R^{n} \rtimes
G(n)$ is an algebraic isomorphism, then $\varphi^{-1}(SO(n,\R))$ is an analytic subgroup of $L$.  The proof is
factored into a sequence of steps.

\begin{lemma}.\\
Let $n \ge 4$.  The centralizer of $\mathscr{M} = \left\{\begin{pmatrix} I_{2} & 0  \\ 0 & M \end{pmatrix}
\ | \ I_{2} = \begin{pmatrix} 1 & 0 \\ 0 & 1 \end{pmatrix} \text{, } M \in SL(n - 2,\R) \right\}$ in
$SL(n,\reals)$ is\\
$\mathscr{A} = \left\{ \begin{pmatrix} A & 0 \\ 0 & D \end{pmatrix} \ | \ A \in GL(2,\reals) \text{, } D = \begin{pmatrix}
\lambda & \  & 0 \\ \  & \ddots & \  \\ 0 & \  & \lambda \end{pmatrix} \in GL(n - 2,\reals) \text{, } \lambda \ne 0
\text{, } \det(A) \cdot \det(D) = 1 \right\}$.\\
If $L$ is a Polish group and $\varphi : L \mapsto \R^{n} \rtimes G(n)$ is an algebraic
isomorphism, then $\varphi^{-1}(\mathscr{A})$ is an analytic subgroup of $L$.
\label{lemma_wolf}
\end{lemma}
Proof:\\
\indent
That $\mathscr{A}$ is as is stated follows from Lemma \ref{lemma_swing} and an elementary computation.
Furthermore\\
$\varphi^{-1}(\mathscr{A}) = \varphi^{-1}(\text{the centralizer of } \mathscr{M} \text{ in } SL(n,\reals)) =$
the centralizer of $\varphi^{-1}(\mathscr{M})$ in $\varphi^{-1}(SL(n,\reals))$ is an\\
analytic subgroup of $L$ by Proposition \ref{proposition_dar}.  $\qed$\\

\begin{lemma}.\\
Let $n \ge 4$ and let $\mathscr{A}$ be as in Lemma \ref{lemma_wolf}.  The commutator subgroup of $\mathscr{A}$
is $\mathscr{B} = \left\{ \begin{pmatrix} B & 0 \\ 0 & I_{n - 2} \end{pmatrix} \ | \ B \in SL(2,\R) \right\}$.
If $L$ is a Polish group and $\varphi: L \mapsto \R^{n} \rtimes G(n)$ is an algebraic isomorphism, then
$\varphi^{-1}(\mathscr{B})$ is an analytic subgroup of $L$.
\label{lemma_gar}
\end{lemma}
Proof:\\
\indent
That $[\mathscr{A},\mathscr{A}] \subset \mathscr{B}$ is a simple computation.  Equality follows from
\cite{jacobson_1985a}.  That $\varphi^{-1}(\mathscr{B})$ is an analytic subgroup of $L$ follows from Lemma
\ref{lemma_wolf} and Lemma \ref{lemma_box}.  $\qed$

\begin{lemma}.\\
Let $\mathbb{F}$ be any field. The centralizer of $\begin{pmatrix} 0 & -1 \\ 1 & 0 \end{pmatrix}$ in
$SL(2,\mathbb{F})$ is $\left\{ \begin{pmatrix} a & -c \\ c & a \end{pmatrix} \ | \ a \text{, } b \in
\mathbb{F} \text{, } a^{2} + c^{2} = 1 \right\}$.
\label{lemma_2dr}
\end{lemma}
Proof:\\
\indent
This follows from a simple computation.  $\qed$

\begin{corollary}.\\
The rotation group $SO(2,\R) = \left\{ \begin{pmatrix} \cos(\theta) & -\sin(\theta) \\ \sin(\theta) &
\cos(\theta) \end{pmatrix} \ | \ \theta \in \R \right\}$ is maximal abelian in $SL(2,\R)$.  In fact $SO(2,\R)$
is the centralizer of $\begin{pmatrix} 0 & -1 \\ 1 & 0 \end{pmatrix}$ in $SL(2,\R)$.
\label{corollary_bar}
\end{corollary}
Proof:\\
\indent
$SO(2,\R)$ is abelian and $\begin{pmatrix} 0 & -1 \\ 1 & 0 \end{pmatrix} \in SO(2,\R)$. Let $B =
\begin{pmatrix} a & b \\ c & d \end{pmatrix}$ be in the centralizer of $\begin{pmatrix} 0 & -1 \\ 1 & 0
\end{pmatrix}$ in $SL(2,\R)$. Lemma \ref{lemma_2dr} implies that $B = \begin{pmatrix} a & -c \\ c & a
\end{pmatrix}$ with $a^{2} + c^{2} = 1$. Hence $a = \cos(\theta)$ and $c = \sin(\theta)$ for some choice of
$\theta \in \R$.  Therefore $B \in SO(2,\R)$ and $SO(2,\R)$ is a maximal abelian subgroup of $SL(2,\R)$.
$\qed$

\begin{lemma}.\\
If $n \ge 2$, then the group $\mathscr{R} = \left\{ \begin{pmatrix} R & 0 \\ 0 & I_{n - 2} \end{pmatrix} \ | \
R \in SO(2,\R) \right\}$ is maximal abelian in\\
$\mathscr{B} = \left\{ \begin{pmatrix} B & 0 \\ 0 & I_{n - 2} \end{pmatrix} \ | \ B \in SL(2,\R) \right\}$. If
$L$ is a Polish group and $\varphi: L \mapsto \R^{n} \rtimes G(n)$ is an algebraic isomorphism, then
$\varphi^{-1}(\mathscr{R})$ is an analytic subgroup of $L$.
\label{lemma_row}
\end{lemma}
\noindent
Proof:\\
\indent
Corollary \ref{corollary_bar} implies that $\mathscr{R}$ is a maximal abelian subgroup in $\mathscr{B}$ if $n
\ge 2$. If $n \ge 4$ then $\varphi^{-1}(\mathscr{B})$ is an analytic subgroup of $L$ by Lemma \ref{lemma_gar}.
If $n = 2$ then $\varphi^{-1}(\mathscr{B})$ is an analytic subgroup of $L$ by Proposition
\ref{proposition_cow}. If $n = 3$ then $\varphi^{-1}(\mathscr{B})$ is an analytic subgroup of $L$ by Lemma
\ref{lemma_nur}. In every case $\varphi^{-1}(\mathscr{R})$ is maximal abelian in $\varphi^{-1}(\mathscr{B})$,
therefore closed in $\varphi^{-1}(\mathscr{B})$ and hence is an analytic subgroup of $L$ since
$\varphi^{-1}(\mathscr{B})$ is an analytic subgroup of $L$.  $\qed$

\begin{lemma}\ \\
Let $H \subset SO(n,\reals)$ be the group generated by elements of the form $\begin{pmatrix} 1& \  &  & \  & \
& \  & \  \\ \  & \ddots & \  & \  & \  & 0 &\  \\ \  & \  & 1 & \  & \ & \  & \  \\ \  & \  & \  & R(\theta)
& \  & \  & \  \\ \  & 0  & \ & \  & 1 & \  & \  \\ \  & \  & \  & \  & \  & \ddots  &\  \\ \  & \  & \  & \ &
\  & \  & 1 \end{pmatrix}$, where $R(\theta)$ is a $2$-dimensional rotation. Then $H$ acts transitively on
$S^{n - 1}$ (and therefore, of course, $SO(n,\R)$ acts transitively on $S^{n - 1}$).
\label{lemma_ton}
\end{lemma}
Proof:\\
\indent
If $x = (x_{1},x_{2}) \ne (0,0)$ and $\|x\| = 1$, there exists a $2$-dimensional rotation matrix $R_{\theta} =
\begin{pmatrix} \cos(\theta) & -\sin(\theta) \\ \sin(\theta) & \cos(\theta) \end{pmatrix} \in SO(2,\R)$ such
that $R_{\theta}\begin{pmatrix} x_{1} \\ x_{2} \end{pmatrix} = \begin{pmatrix} 1 \\ 0 \end{pmatrix}$ for some
$\theta$.

Next, if $x = (x_{1},x_{2},\cdots,x_{n}) \in \R^{n}$ satisfies $\|x\| = 1$, then there exists $R(\theta_{1})
\in SO(2,\reals)$ such that $\begin{pmatrix} 1& \  &  & \  & \  & \  & \  \\ \  & \ddots & \  & \  & \  & 0 &\
\\ \  & \  & 1 & \  & \ & \  & \  \\ \  & \ & \  & \ddots & \  & \  & \  \\ \  & 0  & \ & \  & 1 & \  & \  \\
\  & \  & \  & \  & \  & 1  &\  \\ \  & \  & \  & \  & \  & \  & R(\theta_{1}) \end{pmatrix} \cdot
\begin{pmatrix} x_{1} \\ x_{2} \\ x_{3} \\ \vdots \\ x_{i} \\ \vdots \\ x_{n} \end{pmatrix} = \begin{pmatrix}
x_{1} \\ x_{2} \\ x_{3} \\ \vdots \\ x_{n-2} \\ x'_{n-1} \\ 0 \end{pmatrix} \in S^{n - 1}$. Furthermore there
exists $R(\theta_{2}) \in SO(2,\reals)$ such that $\begin{pmatrix} 1& \  &  & \  & \  & \  & \  \\ \  & \ddots
& \  & \  & \  & 0 &\  \\ \  & \  & 1 & \  & \ & \  & \  \\ \  & \  & \  & \ddots & \  & \  & \  \\ \  & 0  &
\  & \  & 1 & \  & \  \\ \  & \  & \  & \  & \  & R(\theta_{2})  &\  \\ \  & \  & \  & \  & \  & \  & 1
\end{pmatrix} \cdot \begin{pmatrix} x_{1} \\ x_{2} \\ x_{3} \\ \vdots \\ x_{n-2} \\ x'_{n-1} \\ 0
\end{pmatrix} = \begin{pmatrix} x_{1} \\ x_{2} \\ \vdots \\ x_{n-3} \\ x'_{n-2} \\ 0 \\ 0 \end{pmatrix}$.
Continuing in this way, there exists a $R(\theta_{n-1}) \in SO(2,\reals)$ such that $\begin{pmatrix}
R(\theta_{n-1}) & \  &  & \  & \  & \  & \  \\ \  & \ddots & \  & \  & \  & 0 &\  \\ \  & \  & 1 & \  & \ & \
& \  \\ \  & \  & \  & \ddots & \  & \  & \  \\ \  & 0  & \  & \  & 1 & \  & \  \\ \  & \  & \  & \  & \  &
\ddots  &\  \\ \  & \ & \  & \  & \  & \  & 1 \end{pmatrix} \cdot \begin{pmatrix} x_{1} \\ x'_{2} \\ \vdots \\
0 \\ \vdots \\ 0 \\ 0 \end{pmatrix} = \begin{pmatrix} 1 \\ 0 \\ \vdots \\ 0 \\ \vdots \\ 0 \\ 0 \end{pmatrix}
\in S^{n-1}$.

Therefore $H$ acts transitively on $S^{n - 1}$ (and, therefore, $SO(n, \R)$ acts transitively on $S^{n-1}$).
$\qed$

\begin{corollary}\ \\
Identify $SO(n - 1,\R)$ with elements of $SO(n,\reals)$ of the form $\begin{pmatrix} 1 & 0 \\ 0 & A
\end{pmatrix} \subseteq SO(n,\R)$ where $A \in SO(n - 1,\R)$.  Then $S^{n - 1} \simeq SO(n,\R)/SO(n - 1,\R)$
as $SO(n,\R)$ spaces.
\label{corollary_fox}
\end{corollary}
Proof:\\
\indent
This basic fact is well known and can easily be proved using the fact that $SO(n,\reals)$ acts transitively on
$S^{n - 1}$ and that $SO(n - 1,\reals)$ can be identified with the stability group of the north pole.  $\qed$

The next theorem is surely known, but the following proof might be of some interest.

\begin{theorem}\ \\
$SO(n,\R)$ is algebraically generated by elements of the form
$\begin{pmatrix} 1& \  &  & \  & \  & \  & \  \\ \  & \ddots & \  & \  & \  & 0 &\  \\ \  & \  & 1 & \  & \ & \  & \
\\ \  & \  & \  & R(\theta) & \  & \  & \  \\ \  & 0  & \  & \  & 1 & \  & \  \\ \  & \  & \  & \  & \  & \ddots  &\  \\ \
& \  & \  & \  & \  & \  & 1 \end{pmatrix}$, where $R(\theta) \in SO(2,\reals)$.
\label{theorem_car}
\end{theorem}
Proof:\\
\indent
Let $H$ be the subgroup of $SO(n,\R)$ generated by the elements of the form \\ $\begin{pmatrix} 1& \  &  & \ &
\  & \  & \  \\ \  & \ddots & \  & \  & \  & 0 &\  \\ \  & \  & 1 & \  & \ & \  & \ \\ \  & \  & \  &
R(\theta) & \  & \  & \  \\ \  & 0  & \  & \  & 1 & \  & \  \\ \  & \  & \  & \  & \  & \ddots  &\  \\ \ & \ &
\  & \  & \  & \  & 1 \end{pmatrix}$, where $R(\theta) \in SO(2,\reals)$. If $n = 2$, then $H = SO(2,\R)$.
$SO(n - 1,\R) \subset H$ by induction. $H$ and $SO(n,\R)$ act transitively on $S^{n - 1}$ by Lemma
\ref{lemma_ton}. $S^{n - 1} \simeq SO(n,\R)/SO(n - 1,\R)$ as $SO(n,\R)$ spaces by Corollary
\ref{corollary_fox}.  Therefore $H = SO(n,\R)$ by Lemma \ref{lemma_ruf}.  $\qed$

\begin{corollary}\ \\
If $L$ is a Polish group and $\varphi: L \mapsto \R^{n} \rtimes G(n)$ is an algebraic isomorphism, then $\varphi^{-1}
( SO(n,\R))$ is an analytic subgroup of $L$.
\label{corollary_lar}
\end{corollary}
\noindent
Proof:\\
\indent
Let $H$ be the subgroup of $SO(n,\R)$ generated by the elements of the form
$\begin{pmatrix} 1& \  &  & \  & \  & \  & \  \\ \  & \ddots & \  & \  & \  & 0 &\  \\ \  & \  & 1 & \  & \ & \  & \
\\ \  & \  & \  & R(\theta) & \  & \  & \  \\ \  & 0  & \  & \  & 1 & \  & \  \\ \  & \  & \  & \  & \  & \ddots  &\  \\ \
& \  & \  & \  & \  & \  & 1 \end{pmatrix}$, where $R(\theta) \in SO(2,\reals)$.  $H = SO(n,\R)$ by Theorem \ref{theorem_car}.

Let $\mathscr{R}$ be defined as in Lemma \ref{lemma_row}. $\varphi^{-1}(\mathscr{R})$ is an analytic subgroup
of $L$. Let $R_{ij}$ be the result of swapping rows $i$ and $j$ in $I_{n}$. $\det(R_{ij}) = -1$, $R_{ij}^{2} =
I_{n}$, $R_{ij}^{-1} = R_{ij}$ and $R_{ij} \in O(n,\R) = \left\{ A \in GL(n,\reals) \ | \ AA^{T} = A^{T}A =
I_{n} \right\}$, the orthogonal group. Notice that the product of any two $R_{ij}$'s is an element of
$SO(n,\R)$.  Let $\mathscr{P}$ be the finite set of all products of pairs of the $R_{ij}$'s.  It is an
elementary computation that the generators for $H$ are all contained in the finite set of conjugates of
$\mathscr{R}$ by elements of $\mathscr{P}$.  $\varphi^{-1}$ of any of these conjugates is a conjugate of
$\varphi^{-1}(\mathscr{R})$ and is therefore an analytic subgroup of $L$.   Hence, $\varphi^{-1}(SO(n,\reals))
= \varphi^{-1}(H)$ is an analytic subgroup of $L$ by Lemma \ref{lemma_box}.  $\qed$


\newpage

\section{$\R^{n} \rtimes G(n)$ is an Algebraically Determined Polish Group}

This section is devoted to showing that $\R^{n} \rtimes G(n)$ is an algebraically determined Polish group,
thereby proving Theorem \ref{theorem_main}.  This is done by verifyng the hypotheses of Theorem
\ref{theorem_road_map}.

The following proposition is a slight generalization of Lemma 25 of \cite{atim_kallman_2012a} with a different
proof.

\begin{proposition}.\\
Let $\mathscr{K}$ be real or complex inner product space with $\dim(\mathscr{K}) \ge 2$ if $\mathscr{K}$ is
real or $dim(\mathscr{K}) \ge 1$ if $\mathscr{K}$ is complex. Let $x \in \mathscr{K}$ with $\|x\| \le 2$.
Then there exist $y$, $z \in \mathscr{K}$ with $\|y\| = \|z\| = 1$ and $y + z = x$.
\label{proposition_key}
\end{proposition}
\noindent
Proof:\\
\indent
It suffices to consider the case in which $\mathscr{K}$ is an inner product space over the reals, for if
$\mathscr{K}$ is complex, just restrict the scalars to the reals and the real dimension of $\mathscr{K}$ will
be at least two.

Let $x \in \mathscr{K}$ where $\mathscr{K}$ is real and  $\|x\| \le 2$. Choose a unit vector $v$ orthogonal to
$x$ and take $y = \frac{x + av}{2}$ and $z = \frac{x - av}{2}$, where $a = (4 - \|x\|^{2})^{\frac{1}{2}}$.
Clearly $y + z = x$ and $a^{2} = (4 - \|x\|^{2} )$. Therefore $\|y\|^{2} = \| \frac{x+av}{2} \|^{2} =
\frac{1}{4} ( \|x\|^{2} + a^{2} \|v\|^{2}) = \frac{1}{4} ( \|x\|^{2} + (4 - \|x\|^{2}) ) = 1$. Thus $\|y\| =
1$. Similarly one can get $\| z\| = 1$.  $\qed$

\begin{proposition}.\\
Let $\delta > 0$ and let $\mathscr{C} = \left\{ (x,I_{n}) \ | \ \|x\| \le \delta \text{, } x \in \R^{n}
\right\}$. If $L$ is a Polish group and $\varphi: L \mapsto \R^{n} \rtimes G(n)$ is an algebraic
isomorphism, then $\varphi^{-1}(\mathscr{C})$ is analytic in $\varphi^{-1}(\R^{n})$.
\label{proposition_fat}
\end{proposition}
\noindent
Proof:\\
\indent
The statement of this proposition makes sense because $\varphi^{-1}(\R^{n})$ is a Polish group by Lemma
\ref{lemma_bag}.

Fix $x_{0} \in \R^{n}$ with $\|x_{0}\| = \delta/2$. The mapping $\varphi^{-1}(SO(n,\R)) \mapsto
\varphi^{-1}(\R^{n})$,
\begin{eqnarray*}
\varphi^{-1}((0,A)) & \mapsto & \varphi^{-1}((0,A)) \varphi^{-1}((x_{0},I_{n})) \varphi^{-1}((0,A)^{-1}) \\
                    & = & \varphi^{-1}((A(x_{0}),I_{n})))
\end{eqnarray*}
is continuous on $\varphi^{-1}(SO(n,\R))$. Therefore the range of this mapping, namely
$\varphi^{-1}(\left\{(x,I_{n})  \ | \ \|x\| = \delta/2 \right\})$, is an analytic set in $\varphi^{-1}(\R^{n})$
since $\varphi^{-1}(SO(n,\R))$ is an analytic subgroup of $L$ by Corollary \ref{corollary_lar}.

$\varphi^{-1}(\left\{(x,I_{n})  \ | \ \|x\| = \delta/2 \right\}) \times \varphi^{-1}(\left\{(y,I_{n}) \  | \
\|y\| = \delta/2 \right\}) \subset \varphi^{-1}(\R^{n}) \times \varphi^{-1}(\R^{n})$ is also analytic since
the product of two analytic sets is analytic.  Therefore the mapping\\
\begin{displaymath}
\varphi^{-1}((x,I_{n})) \times \varphi^{-1}((y,I_{n})) \mapsto \varphi^{-1}((x,I_{n})) \cdot \varphi^{-1}((y,I_{n})) =
\varphi^{-1} ((x+y,I_{n}))
\end{displaymath}
\begin{displaymath}
\varphi^{-1} ( \left\{(x,I_{n}) : \ \|x\| = \delta/2 \right\}) \times \varphi^{-1} ( \left\{(y, I) : \ \|y\| =
\delta/2 \right\}) \mapsto \varphi^{-1} (\R^{n})
\end{displaymath}
is continuous and therefore has analytic range. But the range of this mapping is
$\varphi^{-1}(\left\{(z,I_{n}) \ | \ \|z\| \leq \delta \right\})$ by Proposition \ref{proposition_key}.
$\qed$

\begin{proposition}.\\
If $L$ is a Polish group and $\varphi: L \mapsto \R^{n} \rtimes G(n)$ is an algebraic isomorphism, then
$\varphi\mid_{\varphi^{-1}(\R^{n})}$ is measurable with respect to $\mathscr{BP}(\varphi^{-1}(\reals^{n}))$.
\label{proposition_car}
\end{proposition}
\noindent
Proof:\\
\indent
The statement of this proposition makes sense because $\varphi^{-1}(\R^{n})$ is a Polish group by Lemma
\ref{lemma_bag}.

Let $B_{\R^{n}}((0,I_{n}),\delta) = \left\{ (x,I_{n}) \ | \ \|x\| < \delta \right\}$. Then
$\varphi^{-1}(B_{\R^{n}}((0,I_{n}),\delta)) = \cup_{\ell \ge 1} \varphi^{-1}(\left\{ (x,I_{n}) \ | \ \|x\| \le
\delta - \frac{1}{\ell} \right\})$ is analytic in $\varphi^{-1}(\R^{n})$ by Proposition
\ref{proposition_fat}.

Fix $x_{0} \in \R^{n}$. Then $\varphi^{-1}(\left\{ (x,I_{n}) \ | \ \|x - x_{0}\| <\delta \right\}) =
\varphi^{-1}((x_{0},I_{n})) \cdot \varphi^{-1}(\left\{ (x,I_{n}) \ | \ \|x\| < \delta \right\})$ is analytic
since $w \mapsto \varphi^{-1}((x_{0},I_{n})) \cdot w$ is continuous on $\varphi^{-1}(\R^{n})$.

Let $\mathscr{O}$ be open in $\R^{n}$ such that $\mathscr{O} = \bigcup_{\ell \ge 1} B_{\R^{n}}
((x_{\ell},I_{n}),\delta_{\ell}) = \bigcup_{\ell \ge 1}\left\{ (x, I_{n})  \mid \| x - x_{\ell} \| <
\delta_{\ell} \right\}$.  Then $\varphi^{-1} (\mathscr{O}) = \bigcup_{\ell \ge 1} \varphi^{-1}( \left\{ (x,
I_{n})  \ | \ \|x - x_{\ell}\| < \delta_{\ell} \right\} )$ is analytic in $\varphi^{-1}(\R^{n})$ since a union
of a sequence of analytic sets is an analytic set.  Hence $\varphi\mid_{\varphi^{-1}( \R^{n})}$ is measurable
with respect to $\mathscr{BP}(\varphi^{-1}(\reals^{n}))$ since analytic sets have the Baire property.
$\qed$\\

\noindent
Proof of Theorem \ref{theorem_main}:\\
\indent
The natural semidirect product $\R^{n} \rtimes G(n)$ is a Polish group. $\varphi^{-1}(\R^{n})$ is closed in
$L$ by Lemma \ref{lemma_bag}. $\varphi^{-1}(G(n))$ is analytic in $L$ in every case by Proposition
\ref{proposition_cow}, by Proposition \ref{proposition_dar} and by Proposition \ref{proposition_rat}. The
natural $\theta$ in every case is always an injection. $\varphi\mid_{\varphi^{-1}(\R^{n})}$ is measurable with
respect to $\mathscr{BP}(\varphi^{-1}(\reals^{n}))$ by Proposition \ref{proposition_car}. Theorem
\ref{theorem_road_map} now implies that $\R^{n} \rtimes G(n)$ is an algebraically determined Polish group.
$\qed$



\end{document}